\documentclass[11pt,twoside]{article}
\usepackage{latexsym}
\usepackage{amssymb,amsbsy,amsmath,amsfonts,amssymb,amscd}
\usepackage{epsfig, graphicx}
\usepackage{hyperref}

\newcommand{\fer}[1]{(\ref{#1})}

\newcommand{\commentout}[1]{}
\newcommand{\R}{\mathbb{R}}
\newcommand{\N}{\mathbb{N}}

\newcommand {\Chi} {{\bf \raise 2pt \hbox{$\chi$}} }

\newcommand {\p}   {\partial}
\newcommand{\beq}{\begin{equation}}
\newcommand{\eeq}{\end{equation}}
\newcommand{\bea} {\begin{array}{rl}}
\newcommand{\eea} {\end{array}}
\newcommand{\bepa}{\left\{ \begin{array}{l}}
\newcommand{\eepa} {\end{array}\right.}

\newcommand{\cqfd}{\begin{flushright}\vspace*{-3mm}$\Box $\vspace{-2mm}\end{flushright}}

\newcommand{\T}{\mathcal T}
\newcommand{\on}{\mbox{ on }}

\newcommand{\eps}{\varepsilon}
\newtheorem{theorem}{Theorem}[section]
\newtheorem{lemma}[theorem]{Lemma}
\newtheorem{definition}[theorem]{Definition}
\newtheorem{remark}[theorem]{Remark}
\newtheorem{proposition}[theorem]{Proposition}
\newtheorem{corollary}[theorem]{Corollary}

\hypersetup{colorlinks=true}
\newcommand{\qed}{{ \hfill
                       {\unskip\kern 6pt\penalty 500
                       \raise -2pt\hbox{\vrule\vbox to 6pt{\hrule width 6pt
                       \vfill\hrule}\vrule} \par}   }}                                            

\begin{document}
\title{Global well-posedness of a conservative relaxed cross diffusion system.}
 \author{Thomas Lepoutre
 \and Michel Pierre \and Guillaume Rolland\\
}
\maketitle
%


\begin{abstract} We prove global existence in time of solutions to relaxed conservative cross diffusion systems governed by nonlinear operators of the form $u_i\to \partial_tu_i-\Delta(a_i(\tilde{u})u_i)$ where the $u_i, i=1,...,I$ represent $I$ density-functions, $\tilde{u}$ is a spatially regularized form of $(u_1,...,u_I)$ and the nonlinearities $a_i$ are merely assumed to be continuous and bounded from below. Existence of global weak solutions is obtained in any space dimension. Solutions are proved to be regular and unique when the $a_i$ are locally Lipschitz continuous. 
\end{abstract}

\section{Introduction}
Introduced by Shigesada {\em et al.} \cite{SKT}, cross diffusion models try to represent the effect of the interaction between species through motion, and not only as usual through reaction. These models have been studied by Levin \cite{L}, Levin and Segel \cite{LS}, Okubo
\cite{Okubo}, Mimura and Murray \cite{MM}, Mimura and Kawasaki
\cite{MK}, Mimura and Yamaguti \cite{MY}, Andreianov {\em et  al.} \cite{ABR}, Bendahmane and Langlais \cite{BL}  and many other authors: a survey by A. J\"{u}ngel may be found in \cite{AJ} for applications to population dynamics.  In those references, a general system is the following:
\beq\label{eq:usual_cross}
\begin{cases}
\p_t u_1-\Delta [u_1(d_1+d_{11}u_1^p+d_{12}u_2^p)]=r_1(u_1,u_2),\\[0.3cm]
\p_t u_2-\Delta [u_2(d_2+d_{21}u_1^p+d_{22}u_2^p)]=r_2(u_1,u_2),\\[0.3cm]
\p_n [u_1(d_1+d_{11}u_1^p+d_{12}u_2^p)]=\p_n[u_2(d_2+d_{21}u_1^p+d_{22}u_2^p)]=0.
\end{cases}
\eeq
For the system \fer{eq:usual_cross} with $p=1$ and Lotka-Volterra-type reaction, there exists a wide literature, studying specific cases of the system where an additional structure keeps it parabolic  or with cross diffusion pressure only on one of the species (see e.g. Wang \cite{Wang} and the many references therein, especially in the introduction). To our knowledge, the most general result on global weak solutions might be found in Chen and J\"{u}ngel \cite{Chen_Jungel} where the entropy structure of the model is used. For  existence of  classical solutions the reader might consult  \cite{Wang, Li_Zhao} by Wang and  Li-Zhao for instance.
In population dynamics, one of the most interesting features of cross diffusion is its effect on steady states: cross diffusion pressure might help the appearance of nonconstant steady states when the reaction structure does not drive to segregation (see Iida-Mimura-Ninomyia \cite{IMN} for instance). However, in these cases, the pattern formation relies on the reaction term (for instance, the convergence to homogeneous steady states in absence of reaction is proved in \cite{Chen_Jungel}).
 
In  \cite{BLMP}, \cite{Lepoutre_PHD}, the first author and collaborators introduced a relaxation of conservative cross diffusion systems, replacing 
\begin{eqnarray*}
\bepa
\p_t u_i-\Delta[a_i(u)u_i]=0, \text{ on }(0,+\infty)\times \Omega, \ \Omega \subset\R^N, \text{ bounded}, \\
u=(u_1,\dots ,u_I),\\
\p_n [a_i(u)u_i]=0 \on (0,+\infty)\times \partial\Omega,\quad u(0,\cdot)=u^0\text{ given, }
\eepa
\end{eqnarray*}
where $a_i:[0,\infty)^I\to [0,\infty)$, by the following relaxed model:
\beq\label{eq:cross_reg}
\bepa
\p_t u_i-\Delta[a_i(\tilde u)u_i]=0, \text{ on }(0,+\infty)\times\Omega, \\
u=(u_1,\dots u_I),\\
\tilde u_i-\delta_i \Delta \tilde u_i=u_i,\text{ on }(0,+\infty)\times \Omega,\ \qquad \delta_i>0,\\
\p_n u_i=\p_n\tilde u_i=0 \text{  on } (0,+\infty)\times \Omega,\quad u(0,\cdot)=u^0\text{ given. }
\eepa
\eeq
This model was introduced in order to investigate the effect of non classical cross diffusion pressure on the segregative behavior (and $a_i(\cdot)$ is often truly nonlinear). One of the purposes was to drive spatial segregation only through motion. Its effects on the stability of the homogeneous equilibria is investigated in \cite{Lepoutre_PHD,BLMP,Lepoutre_Martinez}. This relaxed version is also relevant in some applications: it takes into account that the intensity of the underlying Brownian motion depends on the density of the population measured with a space scale $\delta_i$ and not exactly at the exact location $x$. It takes therefore into account the fact that a species can react to the presence of another species in a neighborhood.  \\
Models with nonlocal diffusion coefficients can be seen also in \cite{BS} (where the self-diffusion coefficients depend on the total population). Nonlocal reaction terms can also be considered (\cite{BNPR,GVA,NPT} for instance), but the goal of our model is more to create patterns only through motion.\\

A first well-posedness result for the relaxed model  was derived in \cite{Lepoutre_PHD}, \cite{BLMP}  in dimension $N=1,2$ and with some restrictions on the structure of the nonlinearities $a_i$ (basically, the $a_i$ are $C^2$ and have at most a polynomial growth in $u$). In this paper, we prove existence of solutions for this system {\em in any dimension and for general nonlinearities $a_i$}, which are only assumed to be continuous and bounded from below. Weak solutions are obtained in general and they are proved to be strong and unique as soon as the $a_i$ are locally Lipschitz continuous. Some $L^2$-estimates are exploited in the spirit of \cite{Pierre_survey} to prove existence of weak solutions. A main point is that $\tilde{u}$ is uniformly bounded in any dimension for these weak solutions. Next, one has to deal with parabolic operators of the form $u_i\to \partial_tu_i-\Delta \left(a_i(\tilde{u})u_i\right)$: they are not of divergence form, but they are uniformly parabolic since $a_i(\tilde{u})$ is then bounded from above and from below. Using the $C^\alpha$-theory for the duals of these operators, namely $U_i\rightarrow \partial_t U_i-a_i(\tilde u)\Delta U_i$, in the spirit of Krylov-Safonov \cite{Krylov_book}, \cite{DG} (see also the book by Lieberman \cite{Lie}), we prove that $\tilde{u}$ is even H\"{o}lder-continuous. This provides continuous coefficients $a_i(\tilde{u})$ for the above operators, and then, $L^p$-estimates classically follow for the solution. When the $a_i$ are locally Lipschitz continuous, even $\partial_tu_i, \Delta\left(a_i(\tilde{u})u_i\right)$ are proved to be in $L^p$ so that the solution is strong: moreover, weak solutions are then proved to be unique.

Let us fix the notations and state the main result. We assume that $\Omega\subset \R^N$ is a bounded subset with a $C^2$-boundary. The exterior normal derivative operator on $\partial\Omega$ is denoted by $\partial_n$. For all $T>0$, we denote $Q_T=(0,T)\times\Omega,\Sigma_T=(0,T)\times\partial\Omega$. For $\alpha\in (0,1]$, we denote 
$$C^\alpha(Q_T)=\{v\in L^\infty(Q_T);\|v\|^{(\alpha)}_T<+\infty\;\},$$
$$\|v\|^{(\alpha)}_T=\|v\|_{L^\infty(Q_T)}+\sup\left\{\frac{|v(t,x)-v(s,y)|}{[|t-s|+|x-y|^2]^{\frac{\alpha}{2}}},\;(t,x),(s,y)\in Q_T\right\}.$$

We will at least assume that
\begin{eqnarray}\label{aicont}
\forall i=1,...,I ,\;\; a_i:[0,\infty)^I\to [0,\infty){\text \;is \;continuous\; and\,:} \inf_{r\in [0,\infty)^I}a_i(r)\geq \underline{d}>0.
\end{eqnarray}
And we are given $\delta_i\in (0,\infty), \forall i=1,...I$. 

\begin{theorem}\label{main} Assume {\rm (\ref{aicont})} and $u^0=(u^0_1,...,u^0_I)\in L^\infty(\Omega,[0,\infty))^I$. Then, there exists a nonnegative solution $u=(u_1,...,u_I)$ to the following problem:
\beq \label{mainsys}
\bepa
\forall T\in (0,\infty), \forall i=1,...,I,\;\forall p\in [1,\infty),\\
u_i\in L^p(Q_T); \tilde{u}_i\in C^\alpha(Q_T)\cap L^p\left(0,T;W^{2,p}(\Omega)\right)\;for\;some\;\alpha\in (0,1],\\
\int_0^ta_i(\tilde u)u_i\in L^p\left(0,T;W^{2,p}(\Omega)\right),\\
 u_i(t)-\Delta[\int_0^ta_i(\tilde u)u_i]=u_i^0 \text{ in }Q_T,\\
\tilde u_i-\delta_i \Delta \tilde u_i=u_i \text{ in} \;Q_T\\
\p_n \left(\int_0^ta_i(\tilde u)u_i\right)=0=\p_n\tilde u_i\;{\text on}\; \Sigma_T.
\eepa
\eeq
If moreover
\beq\label{lip}
\forall i=1,...,I,\; a_i:[0,\infty)^I\to [0,\infty) {\text\;is \;locally\; Lipschitz \;continuous}
\eeq
then, $\forall i=1,...,I,\;\forall T>0, \forall p\in [1,\infty)$,
$$u_i\in L^\infty(Q_T), \forall \tau\in (0,T), \partial_tu_i, \Delta (a_i(\tilde u)u_i)\in L^p\left((\tau,T)\times\Omega\right)$$
and $\partial_tu_i-\Delta (a_i(\tilde u)u_i)=0, \partial_n(a_i(\tilde{u})u_i)=0$ in a pointwise sense. Finally, under assumption {\rm (\ref{lip})},  solutions of {\rm (\ref{mainsys})} are unique.
\end{theorem}

The paper is organized as follows.

\noindent Section \ref{2} first assumes that the nonlinearities $a_i$ are also bounded from above. We prove existence of a weak solution to the system (\ref{mainsys}) by a standard Leray-Schauder fixed-point argument. The underlying space is an adequate subspace of $L^2(Q_T)$ and the required compactness follows essentially from Lemma \ref{pblin}. 

\noindent Section \ref{3} is devoted to the proof of the $L^\infty$-estimate on $\tilde{u}$. Then, the assumption of the bound from above on the $a_i$ may be dropped.

\noindent Section \ref{4} exploits this $L^\infty$-estimate to prove that the weak solution is actually rather regular, and existence as stated in Theorem \ref{main} follows. The $C^\alpha$-theory for nondivergential parabolic operators is used there. An alternative more elementary proof of the regularity is also given when monotonicity properties hold for the $a_i$ together with locally Lipschitz continuity.

\noindent The uniqueness stated in Theorem \ref{main} is proved in Section \ref{5}. It is based on solving an original dual problem, interesting for itself.

\noindent A short Section \ref{6} indicates without proof a complementary approach which provides a constructive and alternative way of proving existence of a solution and which may be used to compute it numerically.

\section{ Global existence when $a_i$ is bounded}\label{2}

In this section, we first prove existence of {\em weak-solutions} of (\ref{mainsys}) on a given interval $[0,T]$ when, besides (\ref{aicont}), the nonlinearities $a_i$ also satisfy
\beq\label{aibounded}
\exists \,\overline{d}>0,\;\;\forall i=1,...,I,\;\;\sup_{r\in [0,\infty)^I}a_i(r) \leq \overline{d}.
\eeq
\begin{proposition}\label{mainproposition}
Let $T>0$. Assume {\rm (\ref{aicont}), (\ref{aibounded})} and $\forall i=1,...,I, u_i^0\in L^2(\Omega;[0,\infty))$. Then,  there exists a nonnegative solution $u=(u_1,...,u_I)$ to the system
\begin{equation}
\left\{
\begin{array}{lll} \label{mainequations}
\;\forall i=1,...,I,\\
u_i\in L^2(Q_T),\;\; \int_0^ta_i(\tilde u)u_i\in L^2\left(0,T;H^2(\Omega)\right),\\
\tilde u_i \in L^2\left(0,T;H^2(\Omega)\right),\; \tilde u_i-\delta_i\Delta \tilde u_i=u_i \;on\; Q_T,\; \tilde{u}_i\geq 0\\
u_i-\Delta(\int_0^t a_i(\tilde u)u_i)=u_i^0 \on Q_T,\\
\partial_n \tilde u_i=0=\partial_n (\int_0^t a_i(\tilde u)u_i) \on \Sigma_T.
\end {array}
\right.
\end{equation}

\end{proposition}

To prove Proposition \ref{mainproposition}, we will use the classical Leray-Schauder's approach, namely (see e.g. \cite{GT}, Theorem 11.3)
\begin{lemma}[Leray-Schauder] \label{LS} Let $(X,\|\cdot\|_X)$ be a Banach space and $\mathcal T:X\rightarrow X$ a continuous compact mapping. Suppose that
$$\exists M>0,\ \forall \sigma \in [0,1],\;\left[\;u\in X,\ u=\sigma \mathcal T u\;\right]\;\Rightarrow\;\left[\;\|u\|_X\leq M\;\right].$$
Then, there exists $u\in X$ such that $u=\mathcal Tu$.
\end{lemma}

To define the mapping $\mathcal T$, we will use the following lemma.

\begin{lemma}\label{pblin}
Let $T>0$, $w_0\in L^2(\Omega;[0,+\infty))$, $A\in L^\infty(Q_T)$,  $\underline a,\ \overline a \in (0,\infty)$ such that $0<\underline a \leq A \leq \overline a<+\infty$. Then there exists a unique nonnegative solution $w=w(A,w_0)$ to
\begin{equation}\label{approximate_problem}
\left\{
\begin{array}{lll}
w\in L^2(Q_T),\;\;\int_0^t Aw\in L^{2}(0,T;H^2(\Omega)),\\
w-\Delta\left(\int_0^tAw\right)=w_0 \on Q_T,\;\partial_n \left(\int_0^t Aw\right )=0 \mbox{ on }\Sigma_T.\\
\end {array}
\right.
\end{equation}
Moreover, if
$$A^n\in L^\infty(Q_T),\;0< \underline{a}\leq A^n\leq \overline{a}<\infty,\;A^n\to A \;a.e.,\;w_0^n\to w_0 \text{ in } L^2(\Omega),$$
then $w(A^n,w_0^n)$ converges strongly in $L^2(Q_T)$ to $w(A,w_0)$.
\end{lemma}

\noindent{\bf Proof of Lemma \ref{pblin}:} Using convolution, we approximate $A$ by a sequence of smooth functions $(A^n)_{n\in \N}\in C^\infty(\overline{Q_T})$ such that $\underline a \leq A^n \leq \overline a$ and $A^n\rightarrow A$ a.e.. Let also $w_0^n$ be a regular approximation of $w_0$. There exists a classical regular nonnegative solution $w^n$ of (see e.g. \cite{LSU}, Theorem V.7.4, applied to the unknown $A^nw^n$)
\begin{equation}\label{localequationN}
 \partial_tw^n-\Delta(A^nw^n)=0 \;on \;Q_T,\;\partial_n (A^n w^n)=0 \;on \;\Sigma_T, w^n(0,\cdot)=w_0^n.
\end{equation}
Integrating (\ref{localequationN}) in time gives
\beq\label{intequN}
w^n(t)-\Delta \left(\int_0^tA^nw^n\right)= w_0^n\;on\;Q_T, \partial_n\left(\int_0^tA^nw^n\right)=0\;on\;\Sigma_T.
\eeq
We multiply by $A^nw^n$ and use the following identity, valid for $z^n=A^nw^n$:
\beq\label{z}
-\int_\Omega z^n\Delta \int_0^tz^n=\int_\Omega\nabla z^n\nabla \int_0^tz^n=\int_\Omega\frac{1}{2}\partial_t|\nabla \int_0^tz^n|^2.
\eeq
We obtain the following estimate after integration in time 
\beq\label{L2n}
\int_{Q_T}A^n(w^n)^2+\int_\Omega\frac{1}{2}|\nabla\int_0^TA^nw^n|^2=\int_{Q_T}w_0^nA^nw^n.
\eeq
In particular
\beq\label{estimL2}
\underline{a}\int_{Q_T}(w^n)^2\leq \overline{a}\sqrt{T}\left(\int_\Omega(w_0^n)^2\right)^{1/2}\left(\int_{Q_T}(w^n)^2\right)^{1/2}\Rightarrow\underline{a}\|w^n\|_{L^2(Q_T)}\leq \overline{a}\sqrt{T}\|w_0^n\|_{L^2(\Omega)}.
\eeq

Now, up to a subsequence, $w^n$ converges weakly in $L^2(Q_T)$ to some $w$. By the pointwise and uniformly bounded convergence of $A^n$ to $A$, for all  $\psi\in L^2(Q_T)$, $\psi A^n$ converges strongly in $L^2(Q_T)$ to $\psi A$ (using the dominated convergence theorem). Thus, $\int_{Q_T}\psi A^nw^n$ converges to $\int_{Q_T}\psi \,A\, w$. In other words, $z^n=A^nw^n$ also converges weakly in $L^2(Q_T)$ to $z=A\,w$.

 By (\ref{intequN}), $\Delta\int_0^tz^n$ is bounded in $L^2(Q_T)$; since $\int_0^tz^n$ is bounded in $L^2(Q_T)$ as well, this implies that $\int_0^tz^n$ is bounded in $L^2(0,T;H^2(\Omega))$. We now may pass to the weak limit in (\ref{intequN}) to deduce that $w$ is solution of (\ref{approximate_problem}). 
 
{\em For the uniqueness}, let $w$ be the difference of two solutions of (\ref{approximate_problem}) (then $w(0)=0$). We denote $S(t)=\int_0^tAw$. {\em Formally}, the idea is to multiply the equation $w-\Delta S=0$ by $S'=A\,w$. Then, after integration
$$\int_{Q_T}\,A\,w^2=\int_{Q_T}S'\Delta S=-\int_{Q_T}\nabla S'\nabla S=-\int_{Q_T}\frac{1}{2}\partial_t|\nabla S(t)|^2=-\int_\Omega\frac{1}{2}|\nabla S(T)|^2\leq 0.$$
Whence $w\equiv 0$ since $A>0$. 

Since we do not know whether $\nabla S'\in L^2(Q_T)$, we have to justify this computation in an approximate way. 
For $h\in (0,T)$, let us denote 
\begin{equation}\label{Sh1}
\left.
\begin{array}{l}
\forall h\in (0,T), \;\;S_h(t):=\frac{S(t+h)-S(t)}{h}=\frac{1}{h}\int_t^{t+h}(Aw)(s)ds.
\end{array}
\right.
\end{equation} 
Note that
\begin{equation}\label{Sh2}
\left.
\begin{array}{l}
S_h\in L^2\left(0,T-h;H^2(\Omega)\right),\; \|S_h-Aw\|_{L^2(Q_{T-h})}\to 0\;as\;h\to 0.
\end{array}
\right.
\end{equation} 
We have
$$\forall t\in [0,T-h),\; w(t+h)+w(t)-\Delta\left[S(t)+S(t+h)\right]=0.$$
We multiply by $S_h(t)$ and integrate over $\Omega$ to obtain
$$\int_\Omega [w(t+h)+w(t)]S_h(t)=-\int_\Omega \nabla S_h(t)[\nabla S(t+h)+\nabla S(t)]=-\int_\Omega \frac{1}{h}\left\{|\nabla S(t+h)|^2-|\nabla S(t)|^2\right\}.$$
After integration on $[0,T-h]$ and an easy change of variable, we have:
\beq\label{h}
\int_{Q_{T-h}} [w(\cdot+h)+w]S_h=-\frac{1}{h}\int_{(T-h,T)\times\Omega}|\nabla S|^2+\frac{1}{h}\int_{(0,h)\times\Omega}|\nabla S|^2\leq \frac{1}{h}\int_{Q_h}|\nabla S|^2.\eeq
To pass to the limit as $h\to 0$, we use 
$$\int_{Q_h}|\nabla S|^2=\int_{Q_h}-S\,w=\int_\Omega-\int_0^h\left[w(t)\int_0^t(Aw)(\sigma)d\sigma\right]dt\leq \|A\|_{L^\infty(Q_T)}h\int_{Q_h}w^2dt.$$
Now, letting $h$ decrease to $0$ in (\ref{h}) and using that $S_h\to Aw$ in $L^2$ (see(\ref{Sh2})), lead to $\int_{Q_T} 2w\,A\,w\leq 0$, whence $w\equiv 0$.

{\em Let us now prove the continuity result}. Let us first notice that, for any solution of (\ref{approximate_problem}), we have the identity
\beq\label{L2A}
\int_{Q_T}A\,w^2+\int_\Omega\frac{1}{2}|\nabla\int_0^TAw|^2=\int_{Q_T}w_0Aw.
\eeq
This may be justified as we did above for the uniqueness (namely in the case $w_0=0$) by passing to the limit in the following identity where $S(t)=\int_0^tA\,w, S_h(t)=[S(t+h)-S(t)]/h$:
\beq
\int_{Q_{T-h}}[w(\cdot+h)+w]S_h+\nabla S_h\nabla [S(\cdot+h)+S]=2\int_{Q_{T-h}}w_0S_h,
\eeq
\beq
\int_{Q_{T-h}}[w(\cdot+h)+w]S_h+\frac{1}{h}\int_{(T-h,T)\times\Omega}|\nabla S|^2-\frac{1}{h}\int_{(0,h)\times\Omega}|\nabla S|^2=2\int_{Q_{T-h}}w_0S_h.
\eeq
And we pass to the limit as above as $h\to 0$ to obtain (\ref{L2A}) (at least a.e.$T$).

Let $w^n=w(A^n,w_0^n)$. As in the beginning of this proof (see (\ref{L2A}), (\ref{estimL2})), the relation
\beq\label{L2nn}
\int_{Q_T}A^n\,(w^n)^2+\int_\Omega\frac{1}{2}|\nabla\int_0^TA^nw^n|^2=\int_{Q_T}w_0^nA^nw^n.
\eeq
proves that $w^n$ is bounded in $L^2(Q_T)$. From (\ref{intequN}), we deduce that  $\int_0^t A^nw^n$ is bounded in $L^2(0,T;H^2(\Omega))$. A subsequence of $\big(w^n,\Delta \int_0^tA^nw^n\big)$ converges weakly in $L^2(Q_T)^2$ to $\big(w,\Delta\int_0^tA\,w\big)$ and $w$ is solution of the limit problem (\ref{approximate_problem}). By uniqueness, the full sequence converges. Since $A^n\to A$ a.e., $\sqrt{A^n}w^n$ converges also weakly in $L^2(Q_T)$ to $\sqrt{A}\,w$ and, by the estimate (\ref{L2nn}), $\nabla\int_0^TA^nw^n$ converges weakly in $L^2(\Omega)$, the limit being necessarily $\nabla \int_0^TA\,w$. In particular
\beq\label{ineqn}
\int_{Q_T}Aw^2\leq \liminf_{n\to\infty}\int_{Q_T}A^n(w^n)^2,\; \int_\Omega|\nabla\int_0^TA\,w|^2\leq \liminf_{n\to\infty}\int_\Omega|\nabla\int_0^TA^nw^n|^2.
\eeq
But, since $\lim_{n\to\infty}\int_{Q_T}w_0^nA^nw^n=\int_{Q_T}w_0A \,w$, and since the identity (\ref{L2A}) is true for $w$, it follows from (\ref{L2nn}), (\ref{L2A}) that {\em equality holds} in the two inequalities (\ref{ineqn}). In particular, the norm of $\sqrt{A^n}w^n$ in $L^2(Q_T)$ converges to the norm of $\sqrt{A}w$; this implies that the $L^2(Q_T)$-weak convergence of $\sqrt{A^n}w^n$ to $\sqrt{A}w$ {\em is actually strong}. Using again the pointwise convergence of $A^n$, we deduce that $w^n$ converges strongly in $L^2(Q_T)$ as well.
\cqfd
\begin{remark}
 As a consequence of $(\ref{L2A})$, there is a constant $C=C(\underline{a},\overline{a}, \|w_0\|_{L^2(\Omega)})$ such that for any solution $w$ of $(\ref{approximate_problem})$,
\begin{equation}\label{estimationL2}
 \|w\|_{L^2(Q_T)}\leq \sqrt{T}C.
\end{equation}
\end{remark}

The next step is the definition of a compact continuous mapping $\mathcal T$ whose fixed points are solutions of $ (\ref{mainequations})$. We introduce the Hilbert space
\begin{equation}\label{X}
X=\Pi_{1\leq i \leq I}X_i,\; X_i=\{v\in L^2(Q_T): \partial_t(J_{\delta_i}v)\in L^2(Q_T)\},
\end{equation}
where the Hilbert norm $\|\cdot\|_i$ is defined on $X_i$ by
$$\|v\|_i^2:=\|v\|^2_{L^2(Q_T)}+\|\partial_t(J_{\delta_i}v)\|^2_{L^2(Q_T)},$$
and where $J_\delta=(I-\delta\Delta)^{-1}$ is the resolvent of the Laplace operator on $L^2(\Omega)$ with homogeneous Neumann boundary conditions, that is
 $$[f\in L^2(\Omega),\; Z=J_\delta f]\Leftrightarrow [Z\in H^2(\Omega),\;Z-\delta\Delta Z=f,\; \partial_nZ=0\;on\;\partial\Omega].$$
\begin{definition} We fix $u^0\in L^2(\Omega, [0,\infty))^I$. Let $v=(v_1,...,v_I)\in X$ and let $\tilde u =(\tilde u_1,...,\tilde u_I)$ be the solution of (see \cite{Brezis}, Proposition 9.24 and Theorem 9.26): 
$$\forall i=1,...,I,\;\tilde{u}_i\in L^2\big(0,T;H^2(\Omega)\big),\;\tilde u_i -\delta_i\Delta \tilde u_i =v_i \mbox{ on }Q_T,\ \ \partial_n \tilde u_i=0 \; on \;\Sigma_T.$$ 
Next, we define 
$$\mathcal{T}: X \rightarrow X\;by\; \mathcal{T}(v):=u=(u_1,...,u_I),$$
where $u_i$ is the solution $w$ of $(\ref{approximate_problem})$ with $A=a_i([\tilde u]^+), w_0=u_i^0$;  $[\tilde{u}]^+=([\tilde{u}_1]^+,...,[\tilde{u}_I]^+)$ and $[\tilde{u}_i]^+$ is the positive part of $\tilde{u}_i$.
\end{definition}

\begin{proposition}\label{prop2} Assume {\rm (\ref{aicont}), (\ref{aibounded})} and $\forall i=1,...,I, u_i^0\in L^2(\Omega;[0,\infty))$.
 Then the mapping $\mathcal T$ is continuous and compact from $X$ into itself.
\end{proposition}
{\bf Proof of Proposition \ref{prop2}: } First, remark that for $v\in X$, $u={\cal T}(v)\in X$. Indeed, since $u_i$ is solution of $(\ref{approximate_problem})$ with $A=a_i([\tilde u]^+)$ and $w_0=u_i^0$, we may write
$$J_{\delta_i}u_i=J_{\delta_i}\Delta \int_0^t Au_i +J_{\delta_i}u_i^0=\int_0^t \Delta J_{\delta_i}(Au_i)+J_{\delta_i}u_i^0 \Rightarrow
\partial_t(J_{\delta_i}u_i)=\Delta J_{\delta_i}(Au_i)\in L^2(Q_T).$$
 Let $v^n$ be a bounded sequence in $X$. Up to a subsequence, me may assume that $v_i^n$ converges weakly to $v_i$ in $L^2(Q_T)$. Then
$$\tilde{u}_i^n-\delta_i\Delta\tilde{u}_i^n=v_i^n\on Q_T, \ \partial _n \tilde u_i^n=0 \on \Sigma_T\;\; \Rightarrow \partial_t\tilde u_i^n=\partial_t(J_{\delta_i}v_i^n).$$
Thus $\tilde{u}_i^n$ is bounded in $L^2\left(0,T;H^2(\Omega)\right)$ and $\partial_t\tilde{u}_i^n=\partial_t\big(J_{\delta_i}v_i^n\big)=J_{\delta_i}(\partial_tv_i^n)$ is bounded in $L^2(Q_T)$. As a consequence, $\tilde{u}_i^n$ is relatively compact in $L^2(Q_T)$, and so is $[\tilde{u}_i^n]^+$. Up to a subsequence again, we may assume that they converge strongly in $L^2(Q_T)$ and a.e. in $Q_T$. By continuity of $a_i$, $a_i([\tilde{u}^n]^+)$ converges a.e. and $0<\underline{d}\leq a_i([\tilde{u}^n]^+)\leq\overline{d}<\infty$. By Lemma \ref{pblin}, $u^n:={\cal T}(v^n)$ converges (up to a subsequence) strongly in $L^2(Q_T)$. Moreover
$$u_i^n=\Delta\left(\int_0^ta_i\left([\tilde{u}^n]^+\right)u_i^n\right)+u_i^0\;\Rightarrow\;\partial_t(J_{\delta_i}u_i^n)=\Delta J_{\delta_i}\left[\left(a_i[\tilde{u}^n]^+\right)u_i^n\right].$$
But the Yosida approximation $\Delta J_{\delta_i}$ is Lipschitz continuous on $L^2(Q_T)$, and $a_i([\tilde{u}^n]^+)u_i^n$ converges in $L^2(Q_T)$. Therefore, $\partial_t(J_{\delta_i}u_i^n)$ converges also in $L^2(Q_T)$. Finally, this proves that $u^n$ converges in $X$ (at least up to a subsequence), whence the compactness of ${\cal T}$.

For the continuity of ${\cal T}$, let $v^n\rightarrow v$ in $X$ as $n\to\infty$.  If $\tilde u^n=(\tilde u_1^n,...,\tilde u_I^n)$ is the solution of
$$\forall\,i=1,...,I,\;\; \tilde u_i^n-\delta_i\Delta \tilde u_i^n=v_i^{n}\on Q_T,\ \partial_n \tilde u_i^n =0 \;on \; \Sigma_T,$$
then $\tilde u_i^n$ converges in $L^2\left(0,T;H^2(\Omega)\right)$ to the solution $\tilde u_i$ of
$$\tilde u_i-\delta_i\Delta \tilde u_i=v_i\on Q_T,\ \partial_n \tilde u_i =0 \; on \; \Sigma_T.$$
By definition, $u^n=\T(v^n)=(u_1^n,...,u_I^n)$ is the solution of
\begin{equation}\label{localequation1}
\left\{
\begin{array}{lll}
u_i^n\in L^2(Q_T),\;\;\int_0^t a_i([\tilde u^n]^+)u_i^n\in L^{2}\left(0,T;H^2(\Omega)\right),\\
u_i^n-\Delta\left(\int_0^ta_i([\tilde u^n]^+)u_i^n\right)=u_i^0 \on Q_T,\;\partial_n \left(\int_0^t a_i([\tilde u^n]^+)u_i^n\right )=0 \mbox{ on }\Sigma_T.\\
\end {array}
\right.
\end{equation}
Using the compactness of $\mathcal T$ proven above, the sequence $(u^n)_{n\in\N}$ is relatively compact in $X$. Let $u^\infty=\lim_{p\to\infty}u^{n_p}$ be a limit point. Up to a subsequence, $\tilde{u}_i^{n_p}$ converges a.e. to $\tilde u_i$. By continuity of $a_i$, $a_i([\tilde u^{n_p}]^+)\rightarrow A_i:=a_i([\tilde u]^+)$ almost everywhere, and it is uniformly bounded from above and from below. According to Lemma \ref{pblin}, we can pass to the limit as $n_p\rightarrow +\infty$ in $(\ref{localequation1})$. By the uniqueness result in Lemma \ref{pblin} with $A=A_i$, we necessarily have $u^\infty=\T(v)$. The sequence $(u^n)_{n\in\N}$ lies in a compact set and has a unique possible limit point, so $u^n=\mathcal T (v^n)\rightarrow \mathcal T(v)$ and $\mathcal T$ is continuous on $X$.
\cqfd

\noindent{\bf Proof of Proposition \ref{mainproposition}:}
Let $T\in (0,\infty)$ and $\sigma \in [0,1]$. Suppose that $u\in X$ is a solution of $u=\sigma \T(u)$. By definition of $\mathcal{T}$, we have
\begin{equation}
\left\{
\begin{array}{lll}
\forall i=1,...,I, u_i\in L^2(Q_T), \;u_i\geq 0,\\
\tilde u_i, \int_0^ta_i(\tilde u)u_i\in L^2(0,T;H^2(\Omega)),\\
\tilde u_i- \delta_i\Delta \tilde u_i=u_i \on Q_T, \partial_n \tilde u_i=0 \on \Sigma_T,\\
u_i-\Delta\int_0^ta_i(\tilde u)u_i=\sigma u_i^0 \on Q_T, \partial _n (\int_0^t a_i(\tilde u)u_i   )=0\on \Sigma_T.
\end {array}
\right.
\end{equation}
The initial conditions $\sigma u_i^0$ are uniformly bounded in $L^2(\Omega)$ for $\sigma\in [0,1]$. Therefore, by the estimate (\ref{estimationL2}),  the function $u_i$ remains bounded in $L^2(Q_T)$, independently of $\sigma$. We also have $\partial_t(J_{\delta_i}u_i)=\Delta J_{\delta_i}(a_i(\tilde u)u_i)$, so $u$ is bounded in $X$ independently of $\sigma$. Using Proposition \ref{prop2} and Leray-Schauder's Lemma \ref{LS},  we can conclude that $\mathcal{T}$ has a fixed point, which is a nonnegative solution of (\ref{mainequations}) (the nonnegativity of $\tilde{u}_i$ is a consequence of $u_i\geq 0$ and of the maximum principle property of $(I-\delta_i\Delta)^{-1}$ with homogeneous Neumann boundary conditions, see e.g. \cite{Brezis}, Proposition 9.30).

\section{$L^\infty$-estimate of $\tilde u$ in Proposition \ref{mainproposition}} \label{3}

A main estimate in the proof of Theorem \ref{main} is given in the next proposition.
\begin{proposition}\label{Linfiniestim}
Assume $u^0\in L^\infty(\Omega,[0,+\infty))^I$ and {\rm (\ref{aicont}), (\ref{aibounded})} as in Proposition \ref{mainproposition}. Let us define
\beq\label{G}
\forall\,k\geq 0,\;\;G(k)=\max_i\{\sup_{r\in [0,k]^I}a_i(r)\}.
\eeq
Then, for any solution $u, \tilde{u}$ of Proposition \ref{mainproposition}, we have
\begin{equation}\label{tildeubound}
\max_{1\leq i\leq I}\left\{\delta_i\|\tilde u_i\|_{L^\infty(Q_T)}+\|\int_0^ta_i(\tilde u)u_i\|_{L^\infty(Q_T)}\right\}\leq M_0+M_1\,T\,G(k_0),
\end{equation}
where $M_0,M_1 \mbox{ and }k_0$ depend only on $u^0,\underline{\delta}:=\min_i\delta_i,\overline{\delta}:=\max_i\delta_i$.
\end{proposition}

The proof of Proposition \ref{Linfiniestim} uses the following classical lemma.

\begin{lemma}\label{Linfinityestimate}
 Let $f\in L^\infty(\Omega)$ and let $w$ satisfy
$$w\in H^2(\Omega),\; w\geq 0,\; \ -\Delta w\leq f \on \Omega,\ \partial_n w=0 \on \partial \Omega.$$
Then there exists $C=C(\Omega)$ such that
\begin{equation}\label{Linftylemma}
\|w\|_{L^\infty(\Omega)}\leq C\left(\|f\|_{L^\infty(\Omega)}+\int_\Omega w\right). 
\end{equation}
\end{lemma}
\noindent{\bf Proof:} First, we rewrite the equation as $w-\Delta w\leq f+w$. Let us fix $p\in (N/2,\infty)$. Using $w\geq 0$, the comparison principle and elliptic regularity theory, we know (see e.g. \cite{GT}, Theorem 8.15) the existence of $C=C(\Omega, p)$ such that
\begin{align*}
\|w\|_{L^\infty}&\leq C  \left(\|f+w\|_{L^p}\right)\leq C\left(\|f\|_{L^p}+\|w\|_{L^p}\right),\\
&\leq C \left(\|f\|_{L^p}+\|w\|_{L^\infty}^{(p-1)/p}(\int_\Omega w)^{1/p}\right),\\
&\leq C \left(\|f\|_{L^p}+\eps\|w\|_{L^\infty}+c(\eps)\int_\Omega w\right) \mbox{ (Young's inequality)}
\end{align*}
and we conclude choosing $\eps$ small enough.
\cqfd
\begin{remark}{\rm Obviously, the conclusion of Lemma \ref{Linfinityestimate} would be the same when assuming only $f\in L^p(\Omega), p>N/2$.
}
\end{remark}
\noindent{\bf Proof of Proposition \ref{Linfiniestim}:} 
We rewrite the equations in $u_i, \tilde u_i$ of  Proposition \ref{mainproposition} as 
\beq\label{eq:integre_temps_reformule}
 \tilde u_i-\Delta\left(\delta_i\tilde u_i +\int_0^t a_i(\tilde u)u_i\right)=u_i^0,\quad \tilde u_i-\Delta w_i=u_i^0,\quad w_i=\delta_i\tilde u_i +\int_0^t a_i(\tilde u)u_i.
\eeq
 We sum up the equations (\ref{eq:integre_temps_reformule}), denoting $\tilde U=\sum_i\tilde u_i, W=\sum_i w_i$\,:
\beq\label{tildeUW}
\tilde U-\Delta W=U^0:=\sum_iu_i^0.
\eeq
Next, we apply Lemma \ref{Linfinityestimate} with $w=W(t), a.e.t, f=U^0$ (note that $-\Delta W(t)\leq U^0$). It gives
\begin{equation}\label{LinfiniW(t)}
a.e.t,\; \|W(t)\|_{L^\infty(\Omega)}\leq C\left(\|U^0\|_{L^\infty(\Omega)}+\int_\Omega W(t)\right). 
\end{equation}
By nonnegativity of $\tilde{u}_i, a_i(\tilde{u})u_i$, we also have (see the definitions of $W,w_i$): $\forall i=1,...,I,\; a.e. t\in [0,T]$:
$$\delta_i\|\tilde u_i(t)\|_{L^\infty(\Omega)}, \|\int_0^ta_i(\tilde u)u_i\|_{L^\infty(\Omega)}\leq \|W(t)\|_{L^\infty(\Omega)}.$$
Then, to end the proof of Proposition \ref{Linfiniestim}, it is sufficient to prove the following lemma.
\cqfd

\begin{lemma}\label{L1estim}
$$a.e. t \in [0,T],\; \int_\Omega W(t)\leq C_0+C_1TG(k_0),$$
where $C_0,C_1,k_0$ depend only on $u^0,\underline{\delta},\overline{\delta}$ and $G$ is defined in {\rm (\ref{G})}.
\end{lemma}

\noindent{\bf Proof of Lemma \ref{L1estim}:}
By integrating the equations on $u_i$ and $\tilde u_i$ in Proposition \ref{mainproposition}, we get:
\beq\label{estimate:L1}
\forall t\geq 0,\qquad \int_\Omega u_i(t)=\int_\Omega \tilde u_i(t)=\int_\Omega u_i^0.
\eeq
 Recall that $\tilde u_i,w_i\in L^2(0,T;H^2(\Omega)), a_i(\tilde u)u_i\in L^2(Q_T)$. We also have $\partial_t\tilde u_i=\Delta J_{\delta_i}(a_i(\tilde{u})u_i) \in L^2(Q_T)$
From (\ref{eq:integre_temps_reformule}), we may write with $\partial_tw_i=\delta_i\partial_t\tilde{u}_i+a_i(\tilde{u})u_i\in L^2(Q_T)$
\begin{equation}\label{wi}
\partial_tw_i-\delta_i\Delta (\partial_tw_i)=a_i(\tilde u)u_i.
\end{equation}
Differentiating $\partial_nw_i=0$ with respect to $t$ on $\partial \Omega$ leads formally to $\partial_n(\partial_tw_i)=0$. Let us check that $\partial_tw_i=\theta(t)$ where $\theta(t)$ is the unique solution of
\begin{equation}\label{thet}
\theta\in L^2\big(0,T;H^2(\Omega)\big),\; a.e. t\in [0,T],\;\theta(t)-\delta_i\Delta \theta(t)=\big(a_i(\tilde{u})u_i\big)(t),\;\partial_n\theta(t)=0\;on\;\partial\Omega.
\end{equation}
Using also $a_i(\tilde{u})u_i\geq 0$, it will then follow that
\begin{equation}\label{wi1}
\partial_tw_i\geq 0,\;\;\partial_tw_i\in L^2(0,T;H^2(\Omega)),\;\;\|\partial_tw_i\|_{L^2(Q_T)}\leq \|a_i(\tilde{u})u_i\|_{L^2(Q_T)}.
\end{equation}
By integration in time of (\ref{thet}), and with $\Theta(t)=\int_0^t\theta(s)ds$, we have
$$\Theta(t)-\delta_i\Delta \Theta(t)=\int_0^t\big(a_i(\tilde{u})u_i\big)(t),\;\partial_n\Theta(t)=0\;on\;\partial\Omega.$$
Comparing with $w_i-\delta_i\Delta w_i=\delta_iu_i^0+\int_0^ta_i(\tilde{u})u_i, \partial_nw_i=0$ implies by uniqueness that:\\ $\Theta(t)=w_i+(I-\delta_i\Delta)^{-1}u_i^0$, whence $\Theta'(t)=\theta=w_i$ after differentiating in $t$.\\

We denote 
$$\widetilde{V}=\sum_i\delta_i\widetilde{u}_i, B=\sum_ia_i(\tilde{u})u_i.$$
 Recall also that 
 $$\widetilde{U}=\sum_i\tilde{u}_i, W=\sum_iw_i, w_i=\delta_i\tilde{u}_i+\int_0^ta_i(\tilde{u})u_i.$$
Summing the $I$ equations in $u_i,\tilde{u}_i$ as in (\ref{tildeUW}), we have
\begin{equation}\label{sumineq}
\overline{\delta}^{-1}\widetilde{V}-\Delta W\leq \widetilde{U}-\Delta W=U^0.
\end{equation}
 We multiply this equation by $\p_t W=\sum_i\partial_tw_i=\p_t\widetilde V+B\geq 0$ (see (\ref{wi1})) and we get 
 $$
 \overline{\delta}^{-1}\int_\Omega \widetilde V(\p_t\widetilde V+B)+\frac{1}{2}\int_\Omega\p_t\left|\nabla W\right|^2\leq \int_\Omega U^0(\p_t\widetilde V+B).
 $$
 We integrate in time to obtain (we denote $\widetilde{V}^0:=\widetilde V(0)=W(0)$)
 \beq\label{eq:sans_nom_inter}
 \int_\Omega \widetilde V^2(T)+\int_{Q_T}2B\widetilde V+\overline{\delta}\int_\Omega |\nabla W(T)|^2\leq \int_\Omega(\widetilde V^0)^2+\overline{\delta}|\nabla \widetilde V^0|^2+2\overline{\delta}U^0(\widetilde V(T)-\widetilde V^0)+\int_{Q_T}2\overline{\delta}U^0B.
 \eeq
 Since we have by definition 
 $$
 \bar\delta U^0=\bar\delta\widetilde U^0-\bar\delta\Delta\widetilde V^0\geq \widetilde V^0-\bar \delta \Delta\widetilde V^0,
 $$
 we have 
 $$
 \int_\Omega (\widetilde V^0 )^2+\bar \delta|\nabla\widetilde V^0 |^2-2\bar\delta U^0\widetilde V^0\leq - \int_\Omega (\widetilde V^0 )^2+\bar \delta|\nabla\widetilde V^0 |^2\leq 0
 $$
 So that \fer{eq:sans_nom_inter} becomes,
  \beq\label{eq:sans_nom}
 \int_\Omega \widetilde V^2(T)+\int_{Q_T}2B\widetilde V+\overline{\delta}\int_\Omega |\nabla W(T)|^2\leq 2\overline{\delta} \int_\Omega U^0\widetilde V(T) +\int_{Q_T}2\overline{\delta}U^0B.
 \eeq
 We have in particular, with $\|U^0\|_\infty=\|U^0\|_{L^\infty(\Omega)}$, and by using (\ref{estimate:L1}):
 \beq\label{eq:1}
 \int_{Q_T} B\widetilde V\leq \overline{\delta}\|U^0\|_\infty\left(\int_\Omega \widetilde V^0+\int_{Q_T}B\right),
 \eeq
Thus, we have for any $k>0$
\beq\label{k}
k\int_{Q_T\cap\{\widetilde V\geq k\}}B\leq \overline{\delta}\|U^0\|_\infty\left(\int_\Omega \widetilde V^0+\int_{Q_T\cap\{\widetilde V < k\}}B+\int_{Q_T\cap\{\widetilde V \geq k\}}B\right).
\eeq
Note that, $\{\widetilde{V}< k\}\subset \cap_i\{\tilde{u}_i\leq k\underline{\delta}^{-1}\}$. Thanks to the $L^1$ estimate (\ref{estimate:L1}), we have  
$$
\int_{Q_T\cap\{\widetilde V < k\}}B=\int_{Q_T\cap\{\widetilde V < k\}}\sum_ia_i(\tilde u)u_i\leq T\left[\int_\Omega U^0\right] G(k\underline{\delta}^{-1}),$$
where $G$ is defined in (\ref{G}).
Finally choosing $k=2\overline{\delta}\|U^0\|_\infty$ in (\ref{k}), we obtain 
$$
\int_{Q_T\cap\{\widetilde V\geq k\}}B\leq \left( 2\int_\Omega \widetilde V^0+ T\left[\int_\Omega U^0\right] G(k_0)\right),\;\; k_0=2\underline{\delta}^{-1}(\overline{\delta}\|U^0\|_\infty).
$$
Adding the two last inequalities gives
\beq\label{eq:bornes_A}
\int_{Q_T} B\leq C_0+C_1T\,G(k_0),
\eeq
where $C_1$ depends only on $u^0, \overline{\delta},\underline{\delta}$.

To end the proof of Lemma \ref{L1estim}, we use that $W(t)=\sum_i\delta_i\tilde{u}_i(t)+\int_0^tB(s)ds$ so that
$$\forall t\in [0,T],\; \int_\Omega W(t)\leq\int_\Omega \widetilde{U}^0+\int_{Q_T}B.$$
\cqfd

From the $L^\infty$-estimate of Proposition \ref{Linfiniestim}, we may now deduce that the problem (\ref{mainsys}) in Theorem \ref{main} has at least a {\em weak solution} under the only assumption of continuity of the $a_i$'s.

\begin{corollary}\label{maincorollary} 
Assume {\rm (\ref{aicont})} (only) and $\forall i=1,...,I, u_i^0\in L^\infty(\Omega;[0,\infty))$. Then,  there exists a nonnegative solution $u=(u_1,...,u_I)$ to the system
\begin{equation}\label{mainequ}
\left\{
\begin{array}{l} 
\forall\, T>0, \;\forall i=1,...,I,\\
u_i, a_i(\tilde{u})u_i\in L^2(Q_T),\;\; \int_0^ta_i(\tilde u)u_i\in L^2\left(0,T;H^2(\Omega)\right),\\
\tilde u_i \in L^\infty(Q_T)\cap L^2\left(0,T;H^2(\Omega)\right),\; \tilde u_i-\delta_i\Delta \tilde u_i=u_i \;on\; Q_T,\; \\
u_i-\Delta(\int_0^t a_i(\tilde u)u_i)=u_i^0 \on Q_T,\\
\partial_n \tilde u_i=0=\partial_n (\int_0^t a_i(\tilde u)u_i) \on \Sigma_T.
\end {array}
\right.
\end{equation}
\end{corollary}

\noindent{\bf Proof:} Here, $a_i$ is assumed to satisfy only (\ref{aicont}) (and not (\ref{aibounded})). Let $T>0$. We introduce
$M_2:=\underline{\delta}^{-1}\left[M_0+M_1TG(k_0)\right]$ where the function $G$ is defined in (\ref{G}) of  Proposition \ref{Linfiniestim}
and $M_0,M_1,k_0$ are defined in (\ref{tildeubound}) of the same proposition. We define
$$\forall r\in [0,M_2]^I,\;\; \overline{a}_i(r):=a_i(r),\forall r\in [0,\infty)^I\setminus [0,M_2]^I, \overline{a}_i(r)=\min\{a_i(r), G(M_2)\}.$$
Then, $\overline{a}_i$ is continuous on $[0,\infty)^I$ and
$$0<\underline{d} \leq \overline{a}_i\leq G(M_2)<\infty,\;\; \overline{a}_i\leq a_i.$$
Therefore, we may apply Proposition \ref{mainproposition} with $a_i$ replaced by $\overline{a}_i$. By Proposition \ref{Linfiniestim}, the corresponding $\tilde{u}$ satisfies 
$$\forall i=1,..,I,\; \|\tilde{u}_i\|_{L^\infty(Q_T)}\leq \underline{\delta}^{-1}\left[M_0+M_1T\overline{G}(k_0)\right],$$
where $\overline{G}$ is defined as in (\ref{G}) with $a_i$ replaced by $\overline{a}_i$. But $\overline{G}(k_0)\leq G(k_0)$, so that
$$\forall i=1,...,I,\; 0\leq \tilde{u}_i\leq M_2,\; \overline{a}_i(\tilde{u})=a_i(\tilde{u}).$$
Therefore, the solution obtained with the data $\overline{a}_i$ is also solution with the data $a_i$.

 This provides a solution of (\ref{mainequ}) in Corollary \ref{maincorollary} with the estimate (\ref{tildeubound}), but only  on $[0,T]$ and it may depend on $T$. To construct a global solution on $(0,\infty)$, we may argue as follows: let $T_p$ be an increasing sequence of times with $\lim_{p\to+\infty}T_p=+\infty$. Let $u^p$ be a solution of our problem on the interval $[0,T_p]$ given by the above proof. For $k\in \N$, we denote by $X^k$ the space $X$ as defined in (\ref{X}) with $T$ replaced by $T_k$ and we denote by $\mathcal T^k:X^k\rightarrow X^k$ the operator $\mathcal T$ with $T=T_k$. For $p\geq k$, we denote $u^{p,k}:=u^p_{[0,T_k]}$ so that $\mathcal T^k(u^{p,k})=u^{p,k}$. We will prove that
 \begin{equation}\label{pk}
 \forall k\in \N,\; (u^{p,k})_{p\geq k}\;{\rm is\;relatively\;compact\;in}\;X^k.
 \end{equation}
 Thus, using a diagonal process, we obtain a sequence $p_m\to\infty$ as $m\to\infty$ and some limit $u$ defined on $(0,\infty)$ so that, for all $k\in \N$, $u^{p_m,k}$ converges to $u_{[0,T_k]}$ in $X^k$ as $m\to\infty$. Then, $\mathcal T_k(u_{[0,T_k]})=u_{[0,T_k]}$ and $u$ is a global solution of  $(\ref{mainequ})$.

  Let $k$ be fixed in $\N$ and let us prove (\ref{pk}). By the $L^\infty$-estimate (\ref{tildeubound}) in Proposition \ref{Linfiniestim}, 
\begin{equation}\label{estimationutildeindependentedep}
\forall\, p\geq k,\;\;\|\tilde{u}_i^p\|_{L^\infty(Q_{T_k})}\leq \frac{1}{\delta_i}[M_0+M_1T_kG(k_0)].
\end{equation}
Thus, $a_i(\tilde{u}^p)$ is uniformly bounded on $Q_{T_k}$. This implies by $(\ref{estimationL2})$ that $u^p$ is  bounded in $L^2(Q_{T_k})^I$ and so is $\partial_t\tilde{u}^p$ since by (\ref{wi1})
$$\delta_i\|\partial_t\tilde{u}_i^p\|_{L^2(Q_{T_k})}\leq 2\|a_i(\tilde{u}^p)u_i^p\|_{L^2(Q_{T_k})}\leq C(k).$$
Thus, $u^{p,k}$ is bounded in $X^k$ and, by compactness of $\mathcal T^k$, it is relatively compact in $X^k$, whence (\ref{pk}).
 \cqfd

\section{Proof of existence in Theorem \ref{main}} \label{4}

Existence of a {\em weak solution} to (\ref{mainsys}) is already proved in Corollary \ref{maincorollary}. It only remains to prove that this solution is actually as regular as stated in Theorem \ref{main}. This will mainly be a consequence of the $L^\infty$-estimate on $\tilde{u}$ proved in the previous section, namely
$$\forall i=1,...,I,\; \|\tilde{u}_i\|_{L^\infty(Q_T)}\leq C_0+C_1T ,\;\; \|a_i(\tilde{u})\|_{L^\infty(Q_T)}\leq C(T),$$
where $C_0,C_1$ depend only on the data and $C(T)=G(C_0+C_1T)$.

We begin by the following simple estimates.
\begin{proposition} \label{estimwi} Let $w_i=\delta_i\tilde{u}_i+\int_0^ta_i(\tilde{u})u_i$ where $u,\tilde{u}$ is solution of (\ref{mainequ}) in Corollary \ref{maincorollary}. Assume $u^0\in L^\infty(\Omega,[0,\infty))^I$. Then, 
\beq\label{wiagain}
\forall T>0,\;\nabla w_i\in L^\infty(Q_T)^N,\; w_i, \partial_tw_i\in L^\infty(Q_T), \;\partial_tw_i\geq 0.
\eeq
\end{proposition}

\noindent{\bf Proof:} The fact that $w_i\in L^\infty(Q_T)$ is a consequence of $(\ref{LinfiniW(t)})$ and Lemma \ref{L1estim}. We recall the two equations (see (\ref{eq:integre_temps_reformule}), (\ref{wi})):
$$ \tilde{u}_i-\Delta w_i=u_i^0,\;\partial_tw_i -\delta_i\Delta(\partial_tw_i)=a_i(\tilde{u})u_i.$$
Since $w_i,\Delta w_i\in L^\infty(Q_T)$ and $\partial_nw_i=0\;on\;\Sigma_T$, we deduce that $\nabla w_i\in L^\infty(Q_T)^N$ (at least).  We have already  seen that $\partial_t w_i\geq 0$ comes directly from the second equation and the nonnegativity of $a_i(\tilde u)u_i$. Now we rewrite this equation as
$$(\partial_tw_i-C(T)\tilde{u}_i)-\delta_i\Delta (\partial_tw_i-C(T)\tilde{u}_i)=(a_i(\tilde{u})-C(T))u_i\leq 0.$$
Together with $\partial_n(\partial_tw_i-C(T)\tilde{u}_i)=0\;on\;\Sigma_T$, this implies
$$\partial_tw_i-C(T)\tilde{u}_i\leq 0, \;so\;that\; 0\leq \partial_tw_i\leq C(T)[C_0+C_1T].$$
\cqfd

We will now prove that $U_i(t,x):=\int_0^t[a_i(\tilde{u})u_i](s,x)ds$ is in $C^\alpha(\overline{Q_T})$ so that, since $\tilde{u}_i=w_i-U_i$, it will follow that $\tilde{u}_i$ is not only bounded, but H\"{o}lder-continuous (at least).

To prove it, we rely on the $C^\alpha$-regularity theory of Krylov-Safonov for the solutions of nondivergence parabolic equations with bounded coefficients. We actually use them in the rather particular case of the operator $-A\Delta$ where $A$ is bounded from above and from below. We may state the result we need as follows:

\begin{lemma}\label{K-S} Let $A\in C(\overline{Q_T}), g\in L^\infty(Q_T), \underline{a},\overline{a}\in (0,\infty)$ with $0<\underline{a}\leq A\leq\overline{a}<\infty$. Let $w\in C^{2,1}(Q_T)\cap C^{1,1}(\overline{Q_T})$ solution of
\beq\label{K-S1}
\left\{
\begin{array}{l}
\partial_tw-A\Delta w=g\;in\;Q_T\\
\partial_nw=0\;on\;\Sigma_T,\; w(0)=0.
\end{array}
\right.
\eeq
Then, there exists $\alpha\in (0,1), C>0$ such that
\beq
\|w\|_T^{(\alpha)}\leq C
\eeq
where $\alpha, C$ depend only on $\underline{a},\overline{a}, T, \|g\|_{L^\infty(Q_T)}, \Omega$.
\end{lemma}
\begin{remark}{\rm Note that the $L^\infty$ estimate of $w$ is easy by a comparison argument (valid here thanks to the a priori regularity of $w$ and of $A$ ): we remark that the function $W(t,x):=t  \,\sup g$ is a supersolution of the problem (\ref{K-S1}), so that $W\geq w$. Doing the same from below, we obtain 
\beq\label{infty}
\|w\|_{L^\infty(Q_T)}\leq T\|g\|_{L^\infty(Q_T)}.
\eeq
Next, we may use the Krylov-Safonov result: the global estimate with homogeneous Neumann boundary conditions as stated above may, for instance, be found in \cite{DG}, Lemma 2.2 (in a quite more general setting).
We more generally refer to \cite{Krylov_book}, \cite{DG}, \cite{Lie} for this kind of results.
}
\end{remark}
We apply this result to prove the regularity of $U_i=\int_0^ta_i(\tilde u)u_i$.
\begin{proposition}\label{Ui} Let $T>0$ and $u^0\in L^\infty(\Omega,[0,\infty))^I$. There exists $\alpha\in (0,1), C>0$ such that
$$\|U_i\|_T^{(\alpha)}+\|\tilde{u}_i\|_T^{(\alpha)}\leq C.$$
\end{proposition}

\noindent{\bf Proof:} Let $u,\tilde{u}$ be the solution of (\ref{mainequ}) in Corollary \ref{maincorollary}. Recall that $0<\underline{d}\leq a_i(\tilde{u})\leq C(T)$. Since Lemma \ref{K-S} a priori applies to regular solutions only, we will use a convenient approximation of $u$. For this, let $A^n$ be a smooth approximation of $a_i(\tilde{u})$  such that
$$0<\underline{d}\leq A^n\leq C(T),\;\; A^n\rightarrow a_i(\tilde{u})\,a.e.$$
Let also $v^n$ be a smooth approximation of $u_i^0$ such that 
$$0\leq v^n\leq \|u_i^0\|_{L^\infty(\Omega)},\;\; v^n\rightarrow u_i^0\;in\;L^2(\Omega).$$
Let $u_i^n$ be the solution of 
$$\partial_tu_i^n-\Delta (A^nu_i^n)=0,\;\partial_nu_i^n=0\;on\;\Sigma_T, u_i^n(0)=v^n.$$
Then, after integration in time, we see that $U_i^n=\int_0^tA^nu_i^n$ satisfies
\beq\label{Uin}
\partial_tU_i^n-A^n\Delta U_i^n=A^nv^n,\;\partial_nU_i^n=0\;on\;\Sigma_T,\; U_i^n(0)=0.\eeq
By Lemma \ref{K-S}, there exists $\alpha,C$ independent of $n$ such that $\|U_i^n\|_T^{(\alpha)}\leq C$. By Lemma \ref{pblin}, $u_i^n$ converges to $u_i$ in $L^2(Q_T)$ which implies that $U_i^n$ also converges to $U_i$ in $L^2(Q_T)$. Whence the estimate of Proposition \ref{Ui} on $U_i$. The estimate on $\tilde{u}_i=w_i-U_i$ follows by combining with Proposition \ref{estimwi} which says that $w_i$ is even Lipschitz continuous.
\cqfd

Now that we know that the coefficient $a_i(\tilde{u})$ is not only bounded but also continuous, we may continue improving the regularity of $u$.

\begin{proposition}\label{Lp} Assume $u^0\in L^\infty(\Omega,[0,\infty))^I$. Then, 
$$\forall p\in [1,\infty),\; \forall T>0,\;\forall i=1,...,I,\; u_i, \partial_tU_i, \Delta U_i\in L^p(Q_T).$$
\end{proposition}

\noindent{\bf Proof:} We may formally write
\beq\label{UVi}
\partial_tU_i-a_i(\tilde{u})\Delta U_i=a_i(\tilde{u})u^0_i,\;\partial_nU_i=0,\; U_i(0)=0.
\eeq
Here $a_i(\tilde{u})$ is continuous on $\overline{Q_T}$ so that, $a_i(\tilde{u})$ being given, this equation has a unique solution : let us call it $V_i$. We set $v_i:=\partial_t V_i/a_i(\tilde{u})$. Then
$$v_i-\Delta V_i=u_i^0, \; V_i=\int_0^ta_i(\tilde{u})v_i,\; \partial_n(\int_0^ta_i(\tilde{u})v_i)=0.$$ Thus, $v_i$ coincides with our $u_i$ (and $V_i$ coincides with our $U_i$) thanks to the uniqueness result of  Lemma \ref{pblin}. 

Moreover, $L^p$-maximal regularity holds for the equation (\ref{UVi}) since $a_i(\tilde{u})$ is continuous (see for instance \cite{LSU}, Theorem 9.1 or \cite{MPS}, Theorem 2.5.2) so that, as $a_i(\tilde{u})u_i^0\in L^\infty(Q_T)\subset L^p(Q_T)$, we have
$$\forall p\in (1,\infty),\; \|\partial_tU_i\|_{L^p(Q_T)},\|\Delta U_i\|_{L^p(Q_T)}\leq C,$$
where $C$ depends on $p, \|a_i(\tilde{u})u^0_i\|_{L^\infty(Q_T)}$ and on the modulus of continuity of the function $a_i(\tilde{u})$. 

Next, from $0 \leq u_i\leq \underline{d}^{-1}a_i(\tilde u)u_i=\underline{d}^{-1}|\partial_tU_i|$, we deduce that $u_i\in L^p(Q_T)$ as well. And $p=1$ is also included since $Q_T$ is bounded.
\cqfd

With Proposition \ref{Lp}, the first part of the existence result in Theorem \ref{main} is now complete. We will now assume that $a_i$ is locally Lipschitz continuous. 

\begin{proposition}\label{fullregul} Besides (\ref{aicont}), assume $a_i$ is locally Lipschitz continuous for all $i=1,...,I$. Assume also $u^0\in L^\infty(\Omega,[0,\infty))^I$. Then
$$\forall i=1,...,I,\;\forall T>0,\;  u_i\in L^\infty(Q_T),\; \forall p<\infty, \forall\tau\in(0,T),\; \partial_t u_i, \Delta(a_i(\tilde{u})u_i)\in L^p((\tau,T)\times\Omega),$$
and
$$\partial_tu_i-\Delta(a_i(\tilde{u})u_i)=0 \on Q_T,\ \partial_n(a_i(\tilde{u})u_i)=0\;on\;\Sigma_T,$$
is satisfied pointwise.
\end{proposition}

\noindent{\bf Proof:} The equation in $u_i$ may also be written (at least formally to start):
\beq\label{aiLip}
\partial_t(a_i(\tilde{u})u_i)-a_i(\tilde{u})\Delta (a_i(\tilde{u})u_i)=u_iDa_i(\tilde{u})\cdot\partial_t\tilde{u}.\\
\eeq
We know that $\delta_i\partial_t\tilde{u}_i+a_i(\tilde{u})u_i\in L^\infty(Q_T)$ (see Proposition \ref{estimwi}), and $a_i(\tilde{u})u_i\in L^p(Q_T)$, for all $p<\infty$, so that $\partial_t\tilde{u}_i\in L^p(Q_T)$ for all $p<\infty$. The right hand side of this equation $F:=u_iDa_i(\tilde{u})\cdot\partial_t\tilde{u}$ is therefore in $L^p(Q_T)$ for all $p<\infty$ since also $Da_i(\tilde{u})\in L^\infty(Q_T)^N$ (as $a_i$ is locally Lipschitz and $\tilde{u}$ is bounded).  

Since $a_i(\tilde{u})$ is continuous on $\overline{Q_T}$, we know (see again e.g. \cite{LSU}, Theorem 9.1 or \cite{MPS}, Theorem 2.5.2), there exists a (unique) solution $\theta$ to
\beq\label{theta}
\left\{
\begin{array}{l}
\forall p<\infty,\; \theta\in C(0,T;L^p(\Omega)), \forall\tau\in(0,T),  \partial_t\theta,\Delta\theta\in L^p((\tau,T)\times\Omega)\\
\partial_t\theta-a_i(\tilde{u})\Delta\theta=F,\;\partial_n\theta=0\;on\;\Sigma_T,\;\theta(0)=a_i(\tilde{u}^0)u_i^0,
\end{array}
\right.
\eeq
and we have
\beq\label{estimtheta}
\|\theta\|_{L^\infty(Q_T)}+\|\partial_t\theta\|_{L^p((\tau,T)\times\Omega)}+\|\Delta \theta\|_{L^p((\tau,T)\times\Omega)}\leq C[\|F\|_{L^p(Q_T)}+\|u_i^0\|_{L^\infty(\Omega)}],
\eeq
where $C$ depends on $\tau,T,p,\Omega$ and of the  modulus of continuity of $a_i(\tilde{u})$.

If we knew that $\theta=a_i(\tilde{u})u_i$, then the proof of Proposition \ref{fullregul} would be complete using moreover:
$$\partial_tu_i=a_i(\tilde{u})^{-1}[\partial_t(a_i(\tilde{u})u_i)-u_iDa_i(\tilde{u})\cdot \partial_t\tilde{u}]\in L^p((\tau,T)\times\Omega).$$
To prove it, we recall (see the proof of Lemma \ref{pblin}) that $u_i$ is the limit of the approximate solutions $u^n$ of
$$\partial_tu^n-\Delta(A^nu^n)=0,\;\partial_nu^n=0\;on\;\Sigma_T,\;u^n(0)=u_i^0,$$
where $A^n$ is smooth and converges pointwise to $a_i(\tilde{u})$ with $0<\min a_i(\tilde{u})\leq A^n\leq \max a_i(\tilde{u})<+\infty$. Moreover, $u^n$ is bounded in $L^p(Q_T)^I$ for all $p<\infty$ by the analysis in Proposition \ref{Lp}. Here, we choose such an approximation $A^n$ which moreover satisfies
$$A^n\to a_i(\tilde{u})\;in\;L^\infty(Q_T),\; \partial_tA^n\to\partial_ta_i(\tilde{u})=Da_i(\tilde{u})\cdot \partial_t\tilde{u}\;in\;L^p(Q_T)\;\forall p<\infty.$$
Then, we apply the estimates (\ref{estimtheta}) to $A^nu^n$ which satisfies
$$\partial_t(A^nu^n)-A^n\Delta(A^nu^n)=u^n\partial_tA^n,\;\partial_nA^n=0\;on\;\Sigma_T,\;A^nu^n(0)=a_i(\tilde{u}^0)u_i^0,$$
and they are preserved at the limit. Whence $\theta=a_i(\tilde{u})u_i$ by uniqueness in (\ref{theta}).
\cqfd

\noindent{\bf Proof of the Existence in Theorem \ref{main}:} It is a consequence of Corollary \ref{maincorollary} and of Propositions \ref{Ui}, \ref{Lp}, \ref{fullregul}.
\cqfd
\begin{remark}\label{regu}{\rm
Note that, not only we proved existence of a solution with the announced regularity, but we even proved that any {\em weak solution} as in Corollary \ref{maincorollary} has actually the announced regularity. This will be useful in the proof of uniqueness}
\end{remark}

\begin{remark}{\rm The assumption that the $a_i$ are bounded from below is essential in our proof of existence, first for the $L^2$-estimate, next to apply the Krylov-Safonov regularity theory. In the case when the $a_i$ degenerate ($a_i\geq 0$), the $L^2$ a priori estimate is to be replaced by  $\sqrt{a_i(\tilde{u})}u_i\in L^2(Q_T)$. However, we loose the $L^2$-compactness of the approximate solutions and also most regularity properties of the solution as well. It could however be interesting to study the possibility of existence of weak solutions.
}
\end{remark}

\begin{remark}{\rm The above analysis relies on the use of the $C^\alpha$-Krylov-Safonov estimates. However, it is interesting to notice that one can prove directly, by an elementary estimate,  that $u\in L^\infty(Q_T)$, without using these estimates in the (rather general situation) where, besides (\ref{aicont}), $a_i$ satisfies
\beq\label{extrai}
\forall i=1,...,I,\; a_i\;{\rm is\;locally\;Lipschitz\;continuous},\; \forall j=1,...,I,\; \partial_{\tilde{u}_j}a_i\geq 0.
\eeq
Once the $L^\infty$-estimate is proved on $u_i$, the full regularity follows by the same arguments as in Proposition \ref{fullregul}.

We indicate below (at least formally) the computations which leads to $u\in L^\infty(Q_T)$.
}
\end{remark}
\noindent{\bf Proof of $u\in L^\infty(Q_T)$ under assumption (\ref{extrai}):}

We write $a_i$ for $a_i(\tilde{u})$ and $a_{ij}=\partial_{\tilde{u}_j}a_i$. We multiply the equation $\partial_tu_i-\Delta(a_iu_i)=0$ by $p(a_iu_i)^{p-1}$ and we integrate over $\Omega$:
\beq\label{firsteq}
\frac{d}{dt}\int_\Omega a_i^{p-1}u_i^p+\int_\Omega p(p-1)(a_iu_i)^{p-2}|\nabla(a_iu_i)|^2=(p-1)\sum_{j}\int_\Omega a_i^{p-2}u_i^pa_{ij}\p_t\tilde u_j.
\eeq
We proved in Proposition \ref{estimwi} that $\partial_t\tilde{u}_j+a_ju_j\leq C(T)<\infty$.
This implies  $\partial_t\tilde{u}_j\leq C(T)$. Plugging this into (\ref{firsteq}), using $a_{ij}\geq 0$, $a_{ij}$ bounded and $a_i\geq\underline{d}$ leads with some $C_T$ independent of $p$ to:
$$
\frac{d}{dt}\int_\Omega a_i^{p-1}u_i^p+\int_\Omega p(p-1)(a_iu_i)^{p-2}|\nabla(a_iu_i)|^2\leq C_T(p-1)\sum_{j}\int_\Omega a_i^{p-1}u_i^p.
$$
Summing over $i$ and using Gronwall's lemma on the term $\sum_i\int_\Omega a_i^{p-1}u_i^p$, we then have  
$$
\sum_i\int_\Omega a_i^{p-1}u_i^p(t)\leq e^{ITC_T(p-1)}\sum_i\int_\Omega a_i^{p-1}u_i^p(0).
$$
Using the lower and upper bounds on $a_i$,  we have with $A, C_T^1$ both independent of $p$:
$$
\sum_i\int_\Omega a_i^{p}u_i^p(t)\leq Ae^{C^1_Tp}(1+\sum_i\|a_iu_i(0)\|_\infty)^p.
$$
This implies
$$
\|(a_iu_i)(t)\|_p\leq A^{1/p} e^{C^1_T}(1+\sum_i\|a_iu_i(0)\|_\infty),
$$
whence the $L^\infty$-estimate on $a_iu_i$ by letting $p\to\infty$, and then on $u_i$ itself by using the lower bound on $a_i$.
\cqfd

\section{Proof of uniqueness in Theorem \ref{main}} \label{5}

Actually, we will prove the more general following result:
\begin{proposition}\label{unique} Let $u^0\in L^\infty(\Omega,[0,\infty))^I$. Assume that for all $i=1,...,I$, $a_i$ satisfies {\rm(\ref{aicont})} and is locally Lipschitz continuous. Then there exists a unique solution to the system {\rm (\ref{mainequ})} in Corollary \ref{maincorollary}.
\end{proposition}

\noindent{\bf Proof:} By Remark \ref{regu}, we already know that any solution of (\ref{mainequ}) satisfies the regularity stated in Proposition \ref{fullregul} and Theorem \ref{main}. Let $u,v$ be two such solutions. We denote $a_i=a_i(\tilde{u}), b_i=a_i(\tilde{v})$. By difference,
$$\partial_t(u_i-v_i)-\Delta\left[ a_i(u_i-v_i)+v_i(a_i-b_i)\right]=0.$$
We set 
$$U_i=u_i-v_i, \tilde{U}_i=\tilde{u}_i-\tilde{v}_i, \widetilde{U}=\tilde{u}-\tilde{v}, A_i=\int_0^1Da_i(t\tilde{u}+(1-t)\tilde{v})dt,$$
so that $a_i-b_i=A_i\cdot(\tilde{u}-\tilde{v})=\sum_jA_{ij}\tilde{U}_j$. Note that $\|A_i\|_{L^\infty}<\infty$. Then
\begin{eqnarray}\label{eqU}
\partial_tU_i-\Delta\left[ a_iU_i+v_iA_i\cdot\widetilde{U}\right]=0,\;\;\partial_\nu( a_iU_i+v_iA_i\cdot\widetilde{U})=0.
\end{eqnarray}

\begin{lemma}\label{dualpb} Let $F\in C_0^\infty(Q_T)^I$. There exists a solution to the dual problem
\beq\label{dualpbeq}
\left\{
\begin{array}{l}
\forall i=1,...,I,\; \varphi_i, \partial_t\varphi_i,\Delta\varphi_i\in L^2(Q_T),\\
\partial_t\varphi_i+a_i\Delta\varphi_i+J_{\delta_i}(B_i\cdot\Delta\varphi)=F_i \on Q_T,\\
\varphi=(\varphi_1,...,\varphi_I), \partial_\nu\varphi_i=0\;on\; \Sigma_T,\;\varphi_i(T)=0,
\end{array}
\right.
\eeq
where $B_i=(B_{i1},...,B_{iI}), B_{ij}=v_jA_{ji}$. 
\end{lemma} 
Assuming this lemma, we multiply each equation (\ref{eqU}) by $\varphi_i$ and we obtain after integration on $Q_T$ (the integrations by parts are allowed, thanks  to the regularity of $u,v,\tilde{u},\tilde{v},\varphi_i$ and the boundary conditions; we also use $\int_{Q_T}U_iJ_{\delta_i}(B_i\cdot\Delta\varphi)=\int_{Q_T}\widetilde{U}_iB_i\cdot\Delta\varphi$):
$$0=\int_{Q_T}U_i[\partial_t\varphi_i+a_i\Delta\varphi_i]+\Delta\varphi_i\,v_iA_i\cdot\tilde{U}=\int_{Q_T}U_iF_i-\tilde{U}_iB_i\cdot\Delta\varphi+\Delta\varphi_i\,v_iA_i\cdot\tilde{U}.$$
Summing these $I$ identities gives $\sum_i\int_{Q_T}U_iF_i=0$ which implies $U\equiv 0$ by arbitrarity of the $F_i$, whence uniqueness.
\cqfd

\noindent{\bf Proof of Lemma \ref{dualpb}:} To solve the dual problem (actually interesting for itself), we may start with $a_i$ replaced by regular approximations $A_i^n$ converging in the usual way to $a_i$ (which means a.e. and uniformly bounded from above and from below), and we first solve
$$\partial_t\theta_i^n+\Delta (A_i^n\theta_i^n)+\Delta J_{\delta_i}(B_i\cdot \theta^n)=\Delta F_i,\;\partial_n(A_i^n\theta_i^n)=0,\; \theta_i(T)=0.$$
This is possible since $\theta\in L^2(Q_T)^I\to \left(\Delta J_{\delta_i}(B_i\cdot\theta)\right)_{1\leq i\leq I}\in L^2(Q_T)^I$
is a Lipschitz perturbation (recall that $B_i\in L^\infty$ and $\Delta J_{\delta_i} $ is the Yosida approximation of the operator $-\Delta$ with homogeneous Neumann boundary conditions). Note that $\int_\Omega\theta_i^n(t)=0$. Next, we solve
$$\Delta \varphi_i^n=\theta_i^n\;in\;\Omega,\; \partial_n(\varphi_i^n)=0\;on\;\partial\Omega, \int_\Omega \phi_i^n=0,$$
so that, "by applying $\Delta^{-1}$" to the equation in $\theta_i^n$, we obtain
\beq\label{eqfinal}
\partial_t\varphi_i^n+A_i^n\Delta\varphi_i^n+J_{\delta_i}(B_i\cdot\Delta\varphi^n)=F_i,\;\partial_n(\varphi_i^n)=0\;on\;\Sigma_T,\; \varphi_i^n(T)=0.
\eeq

Next, multiplying by $\Delta\varphi_i^n$ gives
$$\int_\Omega-\frac{1}{2}\partial_t|\nabla\varphi_i^n|^2+ A_i^n(\Delta\varphi_i^n)^2+\Delta\varphi_i^n J_{\delta_i}(B_i\cdot\Delta\varphi^n) =\int_\Omega F_i\Delta\varphi_i^n\leq \int_\Omega \epsilon (\Delta\varphi_i^n)^2+C_\epsilon F_i^2.
$$
We choose $\epsilon:=\underline{d}/2$ and we deduce
\beq\label{Z}
\int_{\Omega}-\frac{1}{2}\partial_t|\nabla \varphi_i^n|^2+\frac{\underline{d}}{2}(\Delta \varphi_i^n)^2\leq C\int_{\Omega}F_i^2+\int_{\Omega}\nabla Z\nabla \varphi_i^n\leq C\int_{\Omega}F_i^2+\int_\Omega\epsilon|\nabla Z|^2+C_\epsilon|\nabla \varphi_i^n|^2,
\eeq
where $Z-\delta_i\Delta Z= B_i\cdot \Delta\varphi^n,\;\partial_nZ=0$. Multipling this by $Z$ gives
$$\int_\Omega Z^2+\delta_i|\nabla Z|^2=\int_\Omega ZB_i\cdot\Delta\varphi^n\leq \|B_i\|_{L^\infty}\int_\Omega \epsilon Z^2+C_\epsilon|\Delta\varphi^n|^2\Rightarrow\int_\Omega|\nabla Z|^2\leq C\int_\Omega|\Delta \varphi^n|^2.$$
Summing the equations in (\ref{Z}) and choosing adequately $\epsilon$ leads to (with a different $C$)
$$-\partial_t\int_\Omega\sum_i|\nabla\varphi_i^n|^2+\frac{\underline{d}}{2}\sum_i\int_\Omega(\Delta\varphi_i^n)^2\leq C\left[\int_\Omega \sum_i[F_i^2+|\nabla \varphi_i^n|^2]\right].$$
Integrating the Gronwall estimate in $\sum_i|\nabla\varphi_i^n|^2$ and plugging back the terms in $\Delta\varphi_i^n$ yield
$$\sup_{0\leq t\leq T}\int_\Omega\sum_i|\nabla\varphi_i^n|^2+\frac{\underline{d}}{2}\int_{Q_T}|\Delta\varphi^n|^2\leq C\int_{Q_T}|F|^2.$$
By going back to (\ref{eqfinal}), we also obtain that $\partial_t\varphi_i^n$ is bounded in $L^2(Q_T)$. 
Now, we can pass to the limit as $n\to\infty$, weakly in $L^2(Q_T)$ in each term of (\ref{eqfinal}), to prove the existence result of Lemma \ref{dualpb}.
\cqfd
\begin{remark}{\rm We do not know whether uniqueness holds or not without  assuming Lipschitz continuity of the $a_i$. The above proof indicates that uniqueness is essentially equivalent to solving the "dual" problem (\ref{dualpbeq}). The fact that $B_i\in L^\infty(Q_T)$ (which equivalent to the Lipschitz continuity of $a_i$) is strongly used in the estimates to solve (\ref{dualpbeq}). It is not clear how to weaken it.
}
\end{remark}

\section{A constructive approximation procedure} \label{6}

In this Section, we just mention without proof a complementary approach for the proof of existence and regularity of solutions to our system. A interesting point is that it provides a {\em constructive approach } and may be used to provide a numerical approximation schemes  of the solutions (and have actually been used in \cite{Lepoutre_PHD}).

The approximation procedure follows the ideas of \cite{Lepoutre_PHD} to make a time semi-discretization with an explicit treatment of $\tilde u$ and an implicit treatment of $u_i$ in $a_i(\tilde u)u_i$. We fix $T>0$ in the following. Let $\tau>0$ is given. We introduce the following approximate system 
\beq\label{eq:discrete_system}
\left\lbrace\begin{array}{l}
\frac{u_i^{n+1}-u_i^n}{\tau}-\Delta a_i[(\tilde u^n)u_i^{n+1}]=0,\quad\text{in }\Omega,\\[0.3cm]
-\delta_i\Delta \tilde u_i^n+\tilde u_i^n=u_i^n,\quad\text{in }\Omega,\\[0.3cm]
\partial_n \tilde u_i^n=\partial_n u_i^{n+1}=0, \quad\text{on }\partial\Omega.
\end{array}\right.
\eeq
Then, using some adequate discretized version of the approaches in Section \ref{3}, it may be proved that:\\
- there exists a unique solution to the system (\ref{eq:discrete_system}),\\
- the discrete version of the $L^1$-estimate (\ref{eq:bornes_A}) is valid for $B^n=\sum_i a_i(\tilde{u}^{n})u_i^{n+1}$,\\
- an $L^\infty$-estimate may be deduced for $\tilde{u}^n$ ,\\
- the direct $L^\infty$ estimate may be reproduced in this direct situation under the assumption (\ref{extrai}).
Once all these estimates are obtained at the discrete level, one can get Corollary \ref{maincorollary}.

A common issue to both procedures is the need, first to truncate the coefficients (for steps 2 and 3), and then notice \it{a posteriori} that this truncation was useless.

\begin{remark}{\rm
We have used here a situation where the $i$-th relaxation parameter depends on the species (in words $\delta_i$  depends on $i$). It is quite natural on a modeling point of view to consider:
\begin{itemize}
\item
 situations where the relaxation parameter might depend on the species that are observed and on the species that are observing (in words, $a_i=a_i\left(J_{\delta_{i1}}(u_1),\dots ,J_{\delta_{iI}}(u_I)\right)$).
 \item several sensitivities for a single species might have to be considered, short and long distances could matter in the diffusion coefficient (in words two or more different $J_{\delta}(u_j)$ could be taken into consideration in $a_i$).
 \end{itemize}

Our results could be easily extended to these more general situations through rather slight modifications. 
}
\end{remark}

\bibliographystyle{SIAM}

{\tiny
\noindent Thomas Lepoutre\\
INRIA Grenoble Rh\^one-Alpes, team DRACULA\\
Institut Camille Jordan UMR 5208\\
Universit\'e Claude Bernard Lyon 1\\
43 boulevard du 11 novembre 1918\\
69622 Villeurbanne cedex (France)\\
thomas.lepoutre@inria.fr\\

\noindent Michel Pierre\\
ENS Cachan Bretagne, IRMAR, EUB \\
Campus de Ker Lann \\
35170 Bruz (France)\\
michel.pierre@bretagne.ens-cachan.fr\\

\noindent Guillaume Rolland\\
ENS Cachan Bretagne, IRMAR, EUB \\
Campus de Ker Lann \\
35170 Bruz (France)\\
guillaume.rolland@bretagne.ens-cachan.fr

}

\end{document}